\documentclass[12pt]{amsart}

\usepackage{amssymb}

\newcommand{\diag}{\mathop{\rm diag}}

\newtheorem{theorem}{\bf  Theorem}

\begin{document}


\begin{center}
An Elementary Proof of the Polynomial Matrix Spectral
Factorization Theorem
\\[2mm]

 Lasha Ephremidze
          \end{center}

\vskip+0.3cm

 \noindent {\small {\bf Abstract.} A very simple and short proof of the polynomial
 matrix spectral factorization theorem (on the unit circle as well as on the
 real line) is presented, which relies on  elementary complex analysis and linear
  algebra.}

 \vskip+0.1cm\noindent  {\small {\bf Keywords:} Spectral
 factorization, polynomial matrices.
}

\vskip+0.1cm \noindent  {\small {\bf  AMS subject classification
(2010):} 47A68}

\vskip+0.2cm

\section{Introduction}

In this paper, we present an elementary proof of the polynomial
matrix spectral factorization theorem:

\begin{theorem} Let
$$
 S(z)=\sum_{n=-N}^N C_nz^n
$$
be an $m\times m$ matrix function $(${}$C_n\in \mathbb{C}^{m
\times m}$ are matrix coefficients$)$ which is positive definite
almost everywhere on $\mathbb{T}:=\{z\in \mathbb{C}: |z|=1\}$.
Then it admits a factorization
$$
 S(z)=S^+(z)S^-(z),\;\;\;\;\;\;z\in\mathbb{C}\backslash\{0\},
$$
where $S^+(z)= \sum_{n=0}^N A_nz^n$ is an $m\times m$ polynomial
matrix function which is nonsingular inside $\mathbb{T}$,
$\det\,S^+(z)\not=0$ when $|z|<1$, and
$S^-(z)=\overline{S^+\left(1/\overline{z}\right)}^{\,{}_T}=\sum_{n=0}^N
 A^*_nz^{-n}$ is its adjoint, $A_n^*={\overline{A}_n}^{\,{}_T}$,
 $n=0,1,\ldots,N$,
$($respectively, $S^-$ is analytic and nonsingular outside
$\mathbb{T}${}$)$. $S^+$ is unique up to a constant right unitary
multiplier.
\end{theorem}

In the scalar case, $m=1$, the above result is known as the
Fej\'{e}r-Riesz lemma and can be easily proved by considering the
zeroes of $S(z)$.

The matrix spectral factorization $S(t)=S^+(t)\big(S^+(t)\big)^*$,
$|t|=1$, was first established by Wiener [8] in a general case for
any integrable matrix-valued function $S(t)$, $|t|=1$, with an
integrable logarithm of the determinant, $\log \det S(t)\in
L_1({\mathbb T})$. In this case, the spectral factor $S^+$ belongs
to the Hardy space $H_2$.  Wiener proved this theorem by using a
linear prediction theory of multi-dimensional stochastic
processes. A little bit later, by using the same methods,
Rosenblatt [7] showed  that $S^+$ is a polynomial whenever $S$ is
a Laurent polynomial. Since then, many different simplified proofs
of the matrix spectral factorization theorem have appeared in the
literature, see for example [9], [2], [1], [4]. A very short proof
of this result is given in [3], however it uses Wiener's general
spectral factorization existence theorem and some facts from the
theory of the Hardy spaces.

The presented proof  relies only on elementary complex analysis
and linear algebra. This proof is constructive, which makes it
possible to compute the spectral factor approximately at least in
the case of low dimensional matrices (see [5] for a new reliable
computational algorithm of general matrix spectral factorization).
The same pattern can be also used for proving in a straightforward
manner the polynomial matrix spectral factorization theorem on the
real line:

\begin{theorem} Let
$$
 S(z)=\sum_{n=0}^{2N} C_nz^n
$$
be an $m\times m$ matrix function $(${}$C_n\in \mathbb{C}^{m
\times m}$ are matrix coefficients$)$ which is positive definite
almost everywhere on the real line $\mathbb{R}$. Then it admits a
factorization
$$
 S(z)=S^+(z)S^-(z),\;\;\;\;\;\;\;\;z\in\mathbb{C},
$$
where $S^+(z)= \sum_{n=0}^N A_nz^n$ is an $m\times m$ polynomial
matrix function which is nonsingular in the open upper half plane,
$\det\,S^+(z)\not=0$ when ${\rm Im}z>0$, and
$S^-(z)=\overline{S^+\left(\overline{z}\right)}^{\,{}_T}=\sum_{n=0}^N
 A^*_nz^{n}$ is its adjoint, $A_n^*={\overline{A}_n}^{\,{}_T}$,
 $n=0,1,\ldots,N$,
$($respectively, $S^-$ is  nonsingular in the open lower half
plane$)$. $S^+$ is unique up to a constant right unitary
multiplier.
\end{theorem}

Various practical applications of spectral factorization in linear
systems are widely recognized (see, e.g. [6]), where the problem
naturally arises in either of the two different forms commonly
called discrete and continuous. Mathematically both forms are
equivalent under a conformal mapping of the upper half plane into
the unit disk.

\section{Notation}

For $a\in\mathbb{C}$, $a^*=\overline{a}$ denotes its conjugate,
and for a matrix $A$, $A^*=\overline{A}^{\,{}_T}$ denotes its
Hermitian conjugate.

Let $\mathbb{T}_+=\{z\in \mathbb{C}: |z|<1\}$,
$\mathbb{T}_-=\{z\in \mathbb{C}: |z|>1\}\cup\{\infty\}$,
$\mathbb{R}_+=\{z\in \mathbb{C}: {\rm Im} z>0\}$, and
$\mathbb{R}_-=\{z\in \mathbb{C}: {\rm Im} z<0\}$.

$L_1^+(\mathbb{T})$ \big($L_1^-(\mathbb{T})$\big) stands for the
class  of integrable functions  on $\mathbb{T}$ whose Fourier
coefficients with negative (positive) indices equal to zero.

Let $\mathcal{R}^{m\times m}$ be the ring of rational $m\times m$
matrix functions defined in the complex plane. For
$f\in\mathcal{R}^{m\times m} $, the {\em adjoint} matrix function
$\widetilde{f}$ is defined by
$\widetilde{f}(z)=\overline{f\left(1/\overline{z}\right)}^{\,{}_T}$
in the discrete case and by
$\widetilde{f}(z)=\overline{f\left(\overline{z}\right)}^{\,{}_T}$
in the continuous case. Since $f$ is uniquely determined by its
values on $\mathbb{T}$ (on $\mathbb{R}$), and
$\widetilde{f}(z)=(f(z))^*$ for $z\in \mathbb{T}$ (for
$z\in\mathbb{R}$), usual relations for adjoint matrix functions,
like $\widetilde{fg}(z)=\widetilde{g}(z)\widetilde{f}(z)$ and
$\widetilde{f^{-1}}(z)=\widetilde{f}^{-1}(z)$, etc., are valid.
Obviously, if $f(e^{i\theta})\in L_1^+(\mathbb{T})$, then
$\widetilde{f}(e^{i\theta})\in L_1^-(\mathbb{T})$. Whenever
$S^+\in \mathcal{R}^{m\times m}$ is determined, $S^-$ always
denotes its adjoint.

$U\in \mathcal{R}^{m\times m}$ is called paraunitary if
$$
 U(z)\widetilde{U}(z)=I_m
$$
where $I_m$ stands for the $m$-dimensional unit matrix. Note that
$U(z)$ is a usual unitary matrix on the boundary, i.e.
\begin{equation}
U(z)U(z)^*=I_m,\;\;\;z\in \mathbb{T} \;\;\;(z\in \mathbb{R}).
\end{equation}

We say that a matrix function is analytic in a domain if the
entries of the matrix are analytic in the domain.

\section{An elementary proofs of Theorems 1 and 2}

{\em Proof of Theorem $1$}. By Gauss elimination  on the matrix
$S(z)$ and the Fej\'{e}r-Riesz lemma, a factorization
\begin{equation}
S(z)=S_0(z)\widetilde{S_0}(z)
\end{equation}
can be easily achieved with $S_0(z)\in\mathcal{R}^{m\times m}$.
(Namely, if $A=(a_{ij})_{i,j=\overline{1,m}}$ is a positive
definite matrix and $B=(b_{ij})_{i,j=\overline{1,m}}$ is the
unique positive definite matrix such that $A=BB^*$, then the
entries of $B$ can be recursively determined by the formulas
$b_{11}=\sqrt{a_{11}}$, $b_{k1}=a_{k1}/\sqrt{a_{11}}\,^*$,
$k=2,3,\ldots,m$, $b_{nn}=\sqrt{a_{nn}-\sum\nolimits_{j=1}^{n-1}
b_{nj}b_{nj}^*}$, $b_{kn}=\left(a_{kn}-\sum\nolimits_{j=1}^{n-1}
b_{kj}\right)/b_{nn}^*$, $ n=2,3,\ldots,m$, $ k=n+1,\ldots,m$.
These formulas remain valid for rational matrix functions as well.
We need only to assume that $\sqrt{a}$ and $a^*$ are the scalar
spectral factor of $a$ and the adjoint of $a$, respectively.)

If $s_{ij}$ is the $ij$th entry of $S_0$ with a pole at
$a\in\mathbb{T}_+$, then we can multiply $S_0$ by the paraunitary
matrix function $U(z)=\diag[1,\ldots,u(z),\ldots,1]$, where
$u(z)=(z-a)/(1-\overline{a}z)$ is the $jj$th entry of $U(z)$, so
that the $ij$th entry of the product $S_0(z)U(z)$ will not have a
pole at $a$ any longer keeping the factorization (2): $
(S_0U)(z)\,\widetilde{S_0U}(z)=S_0(z)\widetilde{S_0}(z)=S(z). $ In
the same way, one can remove every pole of the entries of $S_0$ at
points in $\mathbb{T}_+$. Thus  $S$ can be represented as a
product
\begin{equation}
 S(z)=S^+_0(z)S_0^-(z),
\end{equation}
where $S^+_0\in\mathcal{R}^{m\times m}$ is  analytic in
$\mathbb{T}_+$

Now, it might happen so that $S_0^+$ is not nonsingular everywhere
on $\mathbb{T}_+$. If $|a|<1$ and $\det S_0^+(a)=0$, then there
exists an $m\times m$ unitary matrix $U$ such that the product
$S_0^+(a)U$ has all 0's in the first column. Hence $a$ is a zero
of every entry of the first column of the matrix function
$S_0^+(z)U$ and the product $S_1^+(z)=
S_0^+(z)U\diag[u(z),1,\ldots,1]$, where
$u(z)=(1-\overline{a}z)/(z-a)$, remains analytic inside
$\mathbb{T}$. While the factorization (3) remains true replacing
$S_0^+$ and $S_0^-$ by $S_1^+$ and $S_1^-$, respectively, the
determinant of $S_1^+$ will have less zeros in $\mathbb{T}_+$ than
the determinant of $S_0^+$. Thus, continuing this process if
necessary, we can remove any singularities in $\mathbb{T}_+$ and
get the factorization
\begin{equation}
 S(z)=S^+(z)S^-(z),
\end{equation}
where $S^+\in\mathcal{R}^{m\times m}$ is  analytic and nonsingular
in $\mathbb{T}_+$.

Now let us show that $S^+$ is in fact a polynomial matrix function
of order $N$. $S^+$ is free of poles on $\mathbb{T}$ since
$S^+(z)\big(S^+(z)\big)^*=S(z)$ for $z\in \mathbb{T}$, and
$z^{-N}S^+(z)=z^{-N}S(z)\big(S^-(z)\big)^{-1}$ is analytic in
$\mathbb{T}_-$. Consequently $S^+$ is analytic in $\mathbb{C}$ and
$z^{-N}S^+(z)$ is analytic in $\mathbb{T}_-$ which implies that
$S^+$ is a polynomial of order $N$.

The proof of the uniqueness of $S^+$ is standard and it is given
only for the sake of completeness. Namely if
$S(z)=S_1^+(z)S_1^-(z)= S_2^+(z)S_2^-(z)$ are two spectral
factorizations of $S(z)$, then $\big(S_2^+(z)\big)^{-1}S_1^+(z)$
is an analytic and nonsingular in $\mathbb{T}_+$ paraunitary
matrix function, which is free of poles and singularities on
$\mathbb{T}$ because of (1). Therefore
$\big(S_2^+(e^{i\theta})\big)^{-1}S_1^+(e^{i\theta})\in
L_1^+(\mathbb{T}) \cap L_1^-(\mathbb{T})$, which implies that it
is a constant matrix function.

\smallskip

{\em Proof of Theorem 2.} This proof can be carried out directly
by the same steps as in the discrete case. The factorization (2)
can be performed exactly in the same way as in the proof of
Theorem 1 assuming under $\sqrt{a}$ a scalar spectral factor of
rational function $a$ which is positive almost everywhere on
$\mathbb{R}$. Elimination of poles of $S_0$ (see (2)) and
singularities of $S_0^+$ (see (3)) in $\mathbb{R}_+$ can be made
by using the paraunitary factors $u(z)=\frac{z-a}{z-\overline{a}}$
and $u(z)=\frac{z-\overline{a}}{z-a}$, respectively. When we have
the factorization (4) where $S^+\in\mathcal{R}^{m\times m}$ is
analytic and nonsingular in $\mathbb{R}_+$, we can prove that
$S^+$ is in fact a polynomial as follows: $S^+$ is free of poles
on $\mathbb{R}$ since $S^+(z)\big(S^+(z)\big)^*=S(z)$ for $z\in
\mathbb{R}$ (as in the discrete case), and $S^+=S(S^-)^{-1}$ is
analytic in $\mathbb{R}_-$. Consequently
$S^+\in\mathcal{R}^{m\times m}$ is analytic in $\mathbb{C}$ and
hence polynomial. The order of $S^+$ is $N$ since if $A$ ia a
nonzero matrix coefficient of the highest order of $z$ in $S^+$,
then $AA^*\not=0$ is the matrix coefficient of the highest order
of $z$ in $S^+S^-$.

The problem of uniqueness of $S^+$ can be reduced to the discrete
case by using the linear fractional transformation
$z\to\frac{i+iz}{1-z}$ which maps $\mathbb{T}_+$ to
$\mathbb{R}_+$.

\end{document}